\documentclass[runningheads]{llncs}
\usepackage{graphicx}
\usepackage{amsmath}
\usepackage{amssymb}

\newcommand{\p}{\partial }

\newcommand{\R}{\ensuremath{\mathbb{R}}}

\newcommand{\mystar}{\mathop{\rm star}}

\begin{document}
\title{Construction of near-boundary Voronoi mesh layers for planar domains\thanks{Supported by RFBR grant 18-01-00726 A}}
%
%\titlerunning{Abbreviated paper title}
% If the paper title is too long for the running head, you can set
% an abbreviated paper title here
%
\author{Vladimir Garanzha\inst{1,2} \and
Liudmila Kudryavtseva\inst{1,2,3} \and Valeriia Tsvetkova\inst{3}
}
\authorrunning{}
% First names are abbreviated in the running head.
% If there are more than two authors, 'et al.' is used.
%
\institute{Dorodnicyn Computing Center FRC CSC RAS, Moscow 119333, Russia  \email{garan@ccas.ru} \url{http://www.ccas.ru/gridgen/lab} \and
Moscow Institute of Physics  and Technology, Moscow, Russia 
\\
 \and
Keldysh Institute of Applied Mathematics RAS, Moscow, Russia\\
\email{liukudr@yandex.ru, lera.tsvetkova@gmail.com}}
\maketitle              % typeset the header of the contribution
\begin{abstract}
We consider problem of constructing purely Voronoi mesh where the union of uncut Voronoi cells approximates the planar computational domain with piecewise-smooth boundary. Smooth boundary fragments are approximated by the Voronoi edges and Voronoi vertices are placed near summits of sharp boundary corners. We suggest self-organization meshing algorithm which covers the boundary of domain by a almost-structured band of non-simplicial Delaunay cells. This band consists of quadrangles on the smooth boundary segment and convex polygons around sharp corners.
Dual Voronoi mesh is double layered orthogonal structure where central line of the layer approximates the boundary. Overall Voronoi mesh has a hybrid structure and consists of high quality convex polygons in the core of the domain and orthogonal layered structure near boundaries. 

\keywords{Voronoi-Delaunay meshing  \and Boundary layer \and Implicit domains.}
\end{abstract}

\section{Introduction}

Construction of hybrid polyhedral meshes in complicated 3d domains is interesting and actively developing field of mesh generation. Well established approach to polyhedral meshing is based on construction of tetrahedral mesh and its approximate dualization \cite{Garimella-2014}, \cite{salome}. In most cases this technique produces high quality polyhedra. Unfortunately near boundary it creates a number of cut cells which should be optimized to get acceptable mesh. Optimality criteria in most cases are contraditory hence costly multicriterial optimization is needed with uncertain outcome.
One can imagine that good solution is construction of Voronoi polyhedral mesh with full uncut Voronoi cells near boundary. We are not aware about such algorithms. Hence, the goal of the paper is to try to construct algorithm which solves above problem in 2d, at least in the practical sense, before treating more complicated 3d case.

Note that approximation of domains by Voronoi tilings and their generalization has rich history, especially in surface reconstruction problems \cite{Amenta-1998}. 
Many algorithms for construction and optimizations of Voronoi meshes were suggested, see \cite{Levy-CVT}, \cite{Alliez-2016}, \cite{Devillers-2010}, \cite {Levy-2010}. Unfortunately these algorithms are not suitable to build Voronoi meshes with regular Voronoi layers near boundaries, which is a topic of present research.

\subsection{Definition of multimaterial implicit domain} 

\begin{figure}
\centerline{(a) \hspace{0.2 \textwidth} (b) \hspace{0.2 \textwidth} (c)} 
\centerline{\includegraphics[width=0.7 \textwidth]{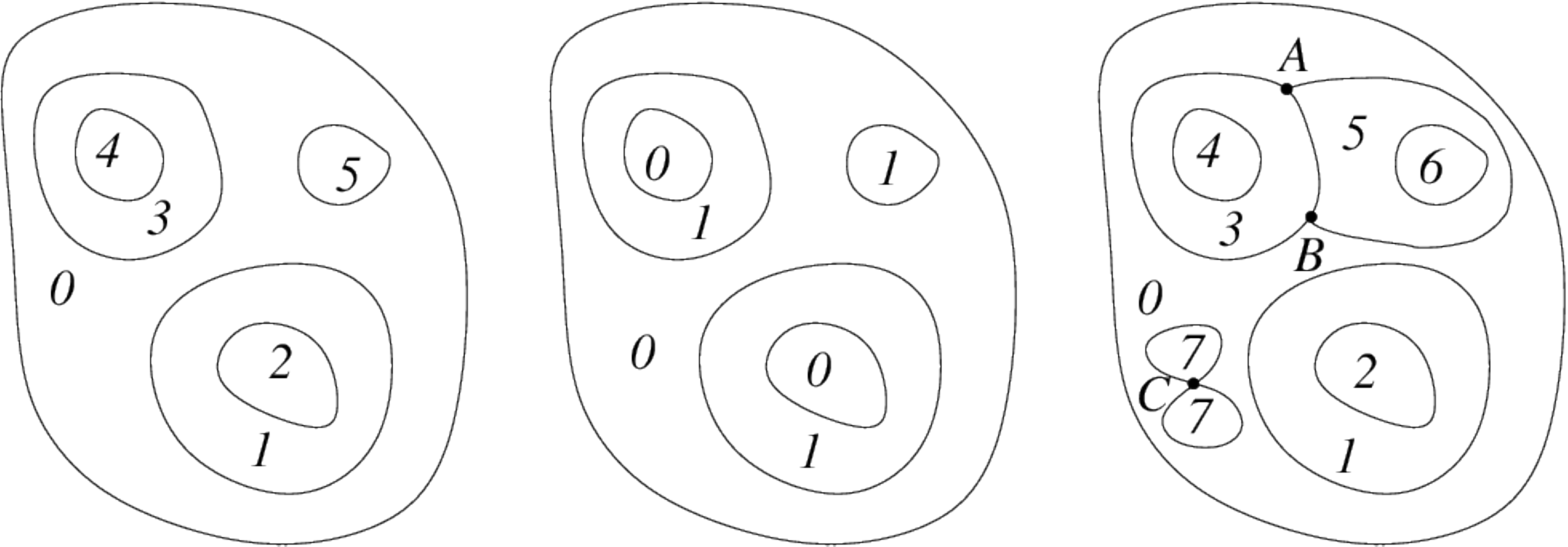}}
\caption{.} \label{fig-multidomain}
\end{figure}

Consider bounded domain $\Omega$ which is partitioned into $N$ subdomains $\Omega_i$, $i = 0, \dots, N-1$. Intuitively one can just imagine a body glued from different materials. We assume that boundary of each subdomain is piece-wise smooth and Lipschitz continuous. The simplest case of multimaterial domain is based on two assumptions: (a) boundary of each subdomain is manifold and (b) multimaterial vertices with neighborhoods containing more than two materials are absent. An example of such a domain is shown in Fig.~\ref{fig-multidomain} (a).

Mesh generation problem in this multimaterial domain is equivalent to mesh generation problem in bimaterial domain shown in Fig.~\ref{fig-multidomain} (a).
One can model such a domain by a single scalar function $u(x): \R^d \to \R$, which is negative inside $\Omega_1$, positive inside $\Omega_0$ and zero isosurface of this function is the boundary. More complicated case is presented in Fig.~\ref{fig-multidomain} (c). Here two multimaterial vertices $A, B$ and non-Lipschitz vertex $C$ are present. Meshing algorithm described below potentially can be applied in this case as well, but we did not tested such a configurations yet.

\begin{figure}
\includegraphics[width=0.8 \textwidth]{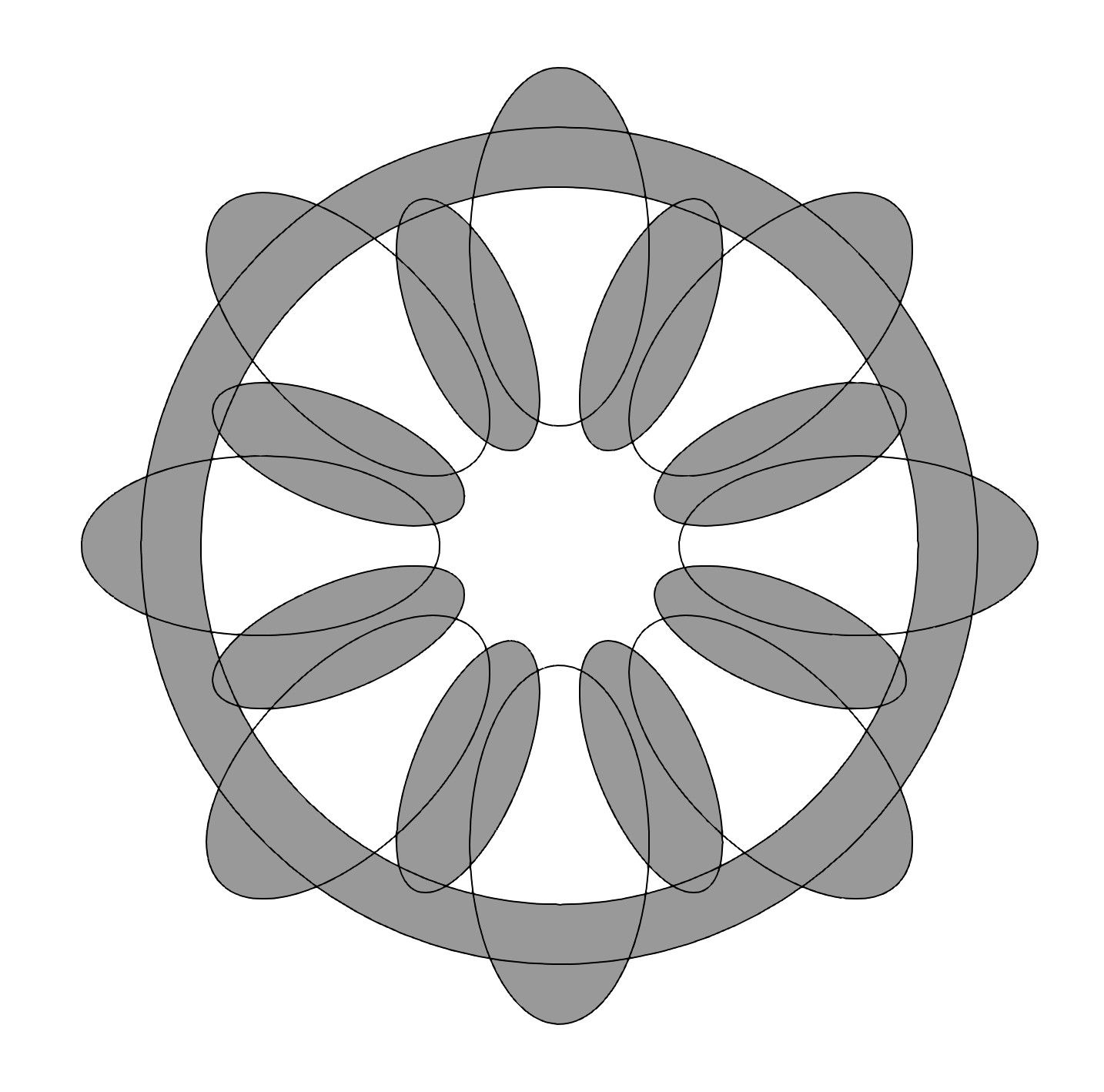}
\caption{Model ``wheel'': construction of implicit domain using Boolean operations.} \label{fig-domain}
\end{figure}

One can use Boolean operations and build quite complicated domains from primitives. Fig.~\ref{fig-domain} shows planar domain that we use as a test case for meshing algorithm.

It is assumed that function $u(x)$ is piecewise smooth, Lipschitz continuous and its derivatives along certain vector field transversal to internal boundary $\Gamma$ are not equal to zero in a finite layer around boundary. In fact it is assumed that behavior of implicit function resembles that of the signed distance function. We can formalize this condition, say, by  assuming the existence of quasi-isometric mapping $y(x): \R^2 \to \R^2$, such that $y (\Omega) = \Omega_y$, and $u(x) = d_s(y(x))$, where $d_s(y)$  is the signed distance function for the surface $y(\Gamma)$. Lipschitz constant for signed distance function cannot exceed unity hence range of values of local Lipschitz constants for function $u(x)$ is defined by the  quasi-
isometry
constants of the mapping $y(x)$.  We do not use this rigorous set of requirements in practice since suggested algorithm is the heuristic one, but always assume that the norm of $\nabla u(x)$, when defined, is bounded from below and from above in a certain layer around $\Gamma$.

\subsection{Voronoi mesh in implicit domain} 

\begin{figure}
\centerline{(a) \hspace{0.4\textwidth} (b)}
\centerline{\includegraphics[width=0.4 \textwidth]{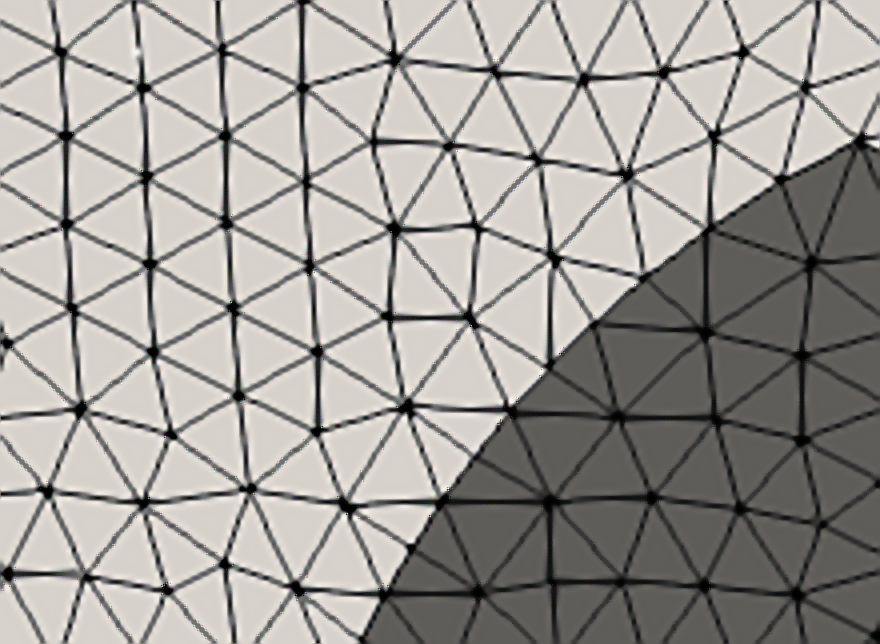}
\includegraphics[width=0.4 \textwidth]{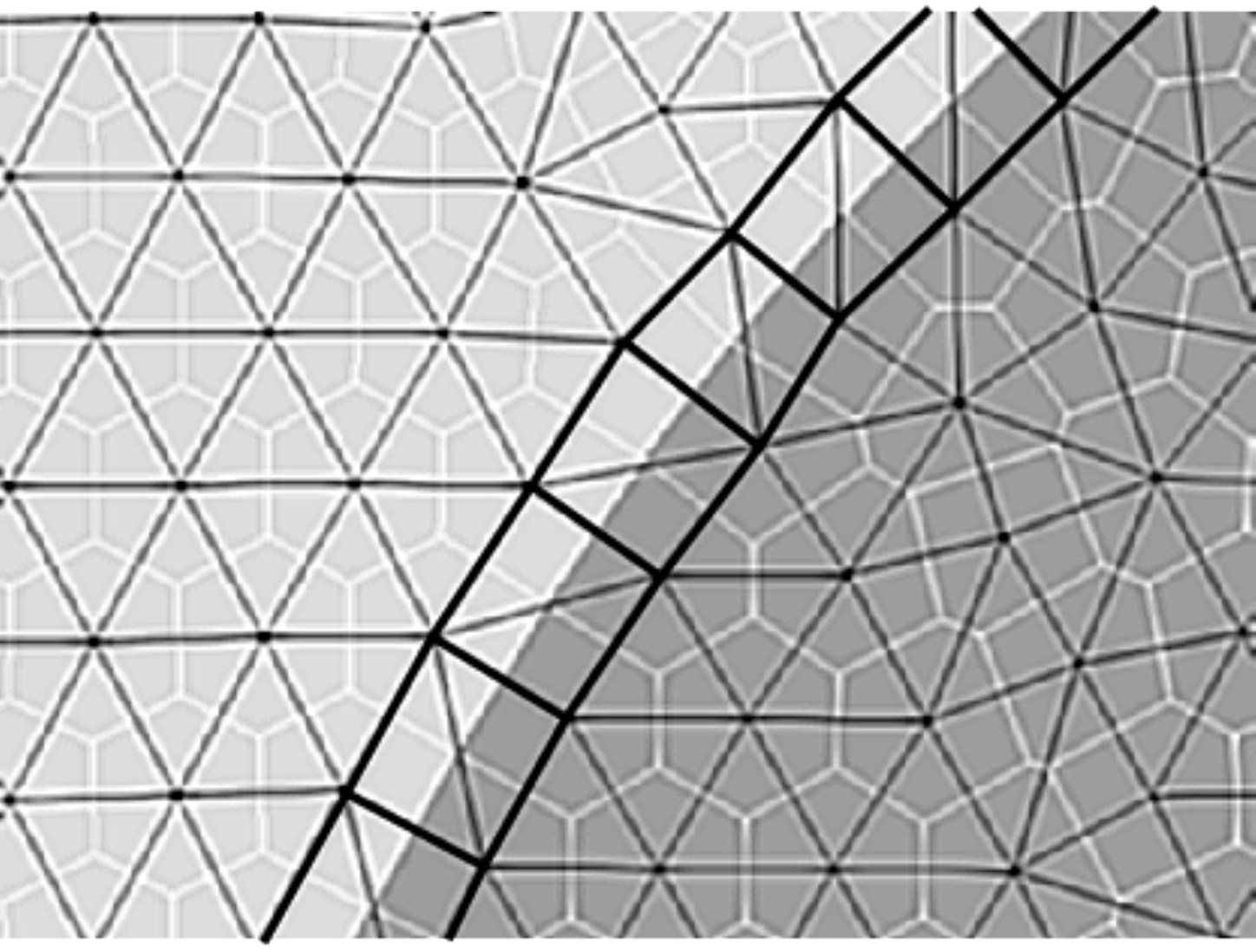}}
\caption{(a) Boundary of domain is approximated by Delaunay edges, (b) boundary of domain is approximated by Voronoi edges.} \label{fig-DV}
\end{figure}

Consider planar mesh $\cal D$ consisting of convex polygons $D_i$ inscribed into circles $B_i$. $D_i$ is convex envelope of all mesh vertices lying on $\p B_i$.
Each circle is empty in a sense that it does not contain any mesh vertices inside. Such a mesh is called Delaunay mesh (Delaunay partitioning). Considering
convex envelope of all centers $c_i$ of circles $B_i$ passing through Delaunay vertex $p_k$ we get Voronoi cell $V_k$. The set of Voronoi cells consitutes what is generally called Voronoi diagram. Since in our setting outer boundary is not approximated, we are not interested in infinite Voronoi cells so we just call resulting object Voronoi mesh. One can approximate internal boundaries using Delaunay mesh as shown in Fig.\ref{fig-DV}(a), or by Voronoi mesh, see Fig.\ref{fig-DV}(b).

Let us briefly explain the difference. Piecewise-smooth boundary $\Gamma$ is approximated by a system of polylines. It is assumed that with mesh refinement polylines converge to $\Gamma$ in the following sense: (a) distance from each straight edge of polyline to certain distinct simple arc of $\Gamma$ should be small; (b) deviation of normal to straight edge from exact normals on the arc should be small; (c) sharp vertices on $\Gamma$ are approximated by sharp vertices on polyline. For Delaunay mesh this polyline is build from Delaunay edges, while for Voronoi mesh polyline is constructed from Voronoi edges. Delaunay edges, dual to the boundary Voronoi edges are orthogonal to the boundary. For smooth fragment of boundary Delaunay cells should be quadrilaterals which make up a band covering the boundary. The midline of this band consisting of Voronoi edges, approximates the boundary as shown in Fig.\ref{fig-DV}(b).
It is well known that all known algorithms generate Delaunay triangulation and not general Delaunay partitions. But edges which split boundary Delaunay cells into triangles has zero dual Voronoi edges and do not influence Voronoi mesh.

\begin{figure}
\centerline{\includegraphics[width=0.5 \textwidth]{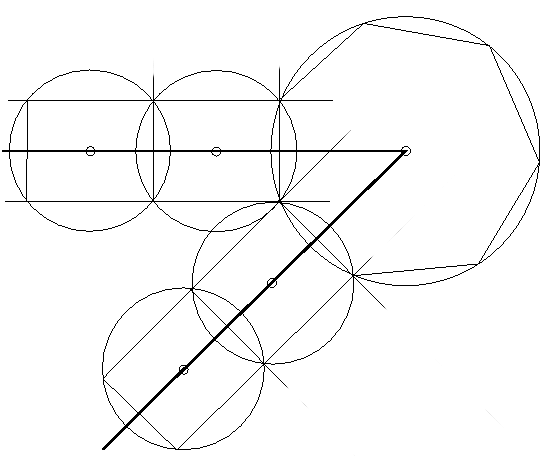}}
\caption{Band of polygonal Delaunay cells and dual Voronoi edges on the boundary of domain.} \label{fig-corner}
\end{figure}

Typical behaviour of Delaunay-Voronoi mesh around sharp boundary vertex is shown in Fig.~\ref{fig-corner}. Regular Delaunay bands consisting of quads are glued together through convex polygonal Delaunay cell. The number of sides in this polygon depends on the sharp vertex angle.

\section{Voronoi meshing algorithm based on self-organization of elastic network} 
In order to build Voronoi meshes in domains with non-smooth boundary we adapt algorithm \cite{GAR-KUDR-2012} which was originaly developed for Delaunay meshing of 2d and 3d implicit domains with piecewise-smooth boundaries. The unknowns in the presented algorithm are  Delaunay mesh vertices which are considered as material points repulsing each other thus modelling elastic medium. Repulsive forces are applied to each pair of vertices belonging to Delaunay edges, i.e. edges with circumferential open balls not containing any other vertices. Each Delaunay edge is treated as compressible strut which tries to expand until prescribed length is reached.   At each step dual Voronoi mesh is constructed and partitioned into two subdomains according to the value of implicit function in the Delaunay vertices. Delaunay mesh is split into three subdomains: subdomain $0$, subdomain $1$ and a set of bands covering the boundary. All Delaunay triangles with circumcenters close to $\Gamma$ are added to the bands. Approximate 
Voronoi boundary polyline is constructed. At this moment mesh refinement is applied provided that local minimum of energy is attained. The idea of mesh refinement is to try to eliminate long Voronoi edges which are not orthogonal to the boundary. It is explained in Fig.~\ref{fig-refine}.

With each Voronoi edge we associate ``sharpening energy'' and ``boundary attraction potential''. Sharpening energy is minimized when Voronoi edge $e$ is orthogonal to the $\nabla u$ at a certain point on $\Gamma$ which we call a target ``touching point'' for $e$. Boundary attraction potential is used as a penalty term for obvious condition that each boundary Voronoi edge is tangential to $\Gamma$ and touches it in certain ``touching point''. We use special variant of preconditioned gradient search method to make one step of minimization. It is convenient to call directions vectors in the minimization technique   ``elastic forces''. 
When due to point displacement under elastic forces edge loses Delaunay property it should be excluded from the list of struts and new Delaunay edges should be created. Hence Voronoi mesh should be rebuilt as well. These steps are repeated until boundary is approximated with reasonable accuracy and correct topology of the near-boundary layers is recovered.

The outcome of the algorithm is certain ''equilibrium`` mesh where elastic forces acting on each point sum to zero. As suggested in \cite{Persson-2004} we build equilibrium mesh in the slightly compressed state.

\subsection{Elastic potential}

Suppose that system of points  \mbox{$\mathcal E = \{p_1, p_2, \dots, p_n \}$} in $\R^2$ is prescribed. Let us denote by $\mathcal T(\mathcal E)$ its Delaunay triangulation.  
We denote by  $\mathcal T_e$ the set of edges of triangulation and by $\mathcal F_b$ the set of Delaunay edges crossing $\Gamma$.  All vertices constitute $2 \times n$ matrix $P$ with $i$-th column equal to $p_i$. We denote the set of near-boundary Delaunay vertices by $P_\Gamma$. Voronoi mesh dual to $\mathcal T$ is denoted by $\mathcal V$, and the set of Voronoi edges detected as a current guess to polyline approximating $\Gamma$ is denoted by $\mathcal E_v$.

With each mesh $\mathcal T$ we associate the following elastic potential
\begin{equation} \label{eq.total-functional}
W(P) =  \theta_r W_r(P) +  \theta_s W_s(P) + \theta_a W_a(P),
\end{equation}
where $W_r(P)$ is the repulsion potential, $W_s(P)$ is the sharpening potential which  serves to align Voronoi boundary edges along isolines of function $u$, $W_a(P)$ is the sharp edge attraction potential.

\paragraph{Repulsion potential.} The repulsion potential is written as follows
\[
W_r(P) = \sum\limits_{e \in \mathcal T_e} w_r(e),
\]
\[
w_r(e) = \left\{ \begin{array}{ll} L_0^2 (\frac{L}{L_0} - 1 -\log (\frac{L}{L_0})) \mbox{ when } L < L_0 \\
0 \mbox{ when } L \geq L_0 \end{array} \right.
\]
where
\[
L = |p_i - p_j|
\]
is the length of the edge  $e$, and $L_0(e)$ is the target length of this edge defined by
\[
L_0(e) = M h(\frac12 (p_i + p_j))
\]
In practice we use $L_0(e) = M \frac12 (h(p_i) + h(p_j))$ in order to diminish number of sizing function calls.

Mesh size distribution is defined by function $f_h(x): \R^2 \to \R$, $f_h(x) > 0$ which can be
interpreted as a relative target edge length at the point $x$. 

\paragraph{Boundary Voronoi edge sharpening and attraction potentials.}
The sharpening functional is written  as follows
\[
W_s(P) = \sum\limits_{e_v \in \mathcal E_v} w_s(e_v),
\]
where the contribution from the boundary Voronoi edge $e_v$ with vertices $c_1, c_2$ looks like
\[
 w_s(e_v) = \frac12 |c_1 - c_2| (n^T (c_2 - c_1))^2
\]
where 
\begin{equation} \label{eq-normal}
n = \frac1{|\nabla u(v^*)|} \nabla u(v^*),
\end{equation}
and $v^*$ is the current approximation of the touching boundary point for the Voronoi edge $e_v$. The simplest choice of $v^*$ is projection of the middle point
\[
c = \frac12 (c_1 + c_2)
\]
of $e_v$ onto $\Gamma$.  

Voronoi edge boundary attraction term is written as
\[
W_a(P) = \sum\limits_{e_v \in \mathcal E_v} w_a(e_v),
\]
where
\[
 w_a(e_v) = \frac12 \left( \frac{L_0}{L} \right)^2 u^2(c)
\]
Here $L$ is the length of the Delaunay edge dual to $e_v$. Hence energy  assigned to shorter Delaunay edges is larger. Since unstable Delaunay edges which  serve to triangulate near-boundary approximate  Delaunay polygons in general are longer compared to stable edges so they produce smaller constribution to total energy and have small influence on positions of vertices.

\subsection{``Elastic forces'' and practical iterative algorithm}

It is convenient to introduce the notions of ``repulsive  forces'',   ``sharpening forces'' and ``boundary attraction forces'' which denote the contribution to the direction vector from the repulsion, sharpening and boundary attractions terms, respectively.

Roughly speaking, these ``forces'' are introduced as follows
\begin{equation} \label{eq.final-forces}
\begin{split}
 \delta p^k_i =  -\frac{\theta_r}{{d_r}^k_i}\frac{\p W_r}{\p p_i}(P^k)  - \frac{\theta_s}{ {d_s}^k_i} \frac{\p W_s}{\p p_i}(P^k) - \frac{\theta_a}{ {d_a}^k_i} \frac{\p W_a}{\p p_i}(P^k) =  \\ 
 = F_e(p^k_i) +  F_s(p^k_i) + F_a(p^k_i),
 \end{split}
\end{equation}
Here $k$ is the iteration number, $p_i$ is the $i$-th vertex in the Delaunay mesh $P^k$,
${d_r}^k_i, {d_s}^k_i, {d_a}^k_i$ are the scaling factors.

Since Newton law is not used to describe the motion of mesh vertices these ``forces'' are speculative and are just used to facilitate intuitive understanding of the algorithm.

In order to present precise formulae for computation of forces it is convenient to introduce the following notations. Let $\mystar_e(p_i)$ denote the set of the mesh edges originating from the vertex $p_i$, while $\mystar(p_i)$ will denote the set of vertices of these edges excluding $p_i$.  In all cases we assume that every boundary star is ordered, i.e. its entities are numbered counterclockwise around $p_i$ looking from outside the domain. Below we omit upper index $k$.

\paragraph{Repulsive ``force''.} For internal vertex $p_i$
\[
{F}_r(p_i) = -\frac{\theta_r}{d_i} {\sum_{p_j \in \mystar{p_i}} \phi_r (p_i, p_j) (p_i - p_j)}, \
{d_r}_i = {\sum_{p_j \in \mystar{p_i}} \phi_r(p_i, p_j)},
\]
where
\[
\phi_r (p_i, p_j) = (\frac{L_0}{L} - 1) \frac{L_0}{L}, \ L = |p_i - p_j|, \ L_0 = Mh(\frac12 (p_i + p_j))
\]

\paragraph{Sharpening and boundary attraction forces.} Sharpening force can be written as follows 
\[ 
{F}_s(p_i) = -\frac {\sum\limits_{e_v: \ p_i \in \mathop{\rm dual} e_v} \Pi_r (q |c_1 - c_2| n^T (c_1 - c_2))}{\sum\limits_{e_v: p_i \in \mathop{\rm dual} e_v} |c_1 - c_2| |q|^2},
\]
here $c_1, c_2$ are vertices of the edge $e_v$,  $c = \frac12 (c_1 + c_2)$,
vector $n$ is defined in (\ref{eq-normal}), and
\[
q = (C_2 - C_1)^T n, \ C_1 = \frac{\p c_1}{\p p_i}, \ C_2 = \frac{\p c_2}{\p p_i}
\]
In order to write down expression for matrix $C_1$, consider Delaunay triangle $T_1$ with counterclockwise ordered vertices $p_i, p_j, p_k$ whose circumcenter is $c_1$. Then
\[
C_1^T = (c_1 - p_i \, c_1 - p_i) (p_j - p_i \, p_k - p_i)^{-1}
\]
Formula for $C_2$ is similar.

Nonlinear operator $\Pi_r$  is responsible for interaction between repulsive force and sharpening force.

Consider contribution to ${F}_s(p_i)$ from Voronoi edge $e_v$. Denote by $e = (p_j - p_i) / |p_j - p_i|$, where $p_i, p_j$ are vertices of Delaunay edge dual to $e_v$.
If 
\[
e^T F_r(p_i) e^T q n^T (c_1 - c_2) < 0 
\]
then
\[
\Pi_r(q n^T (c_1 - c_2)) = q n^T (c_1 - c_2) - e e^T q n^T (c_1 - c_2)
\]
otherwise
\[
\Pi_r(q n^T (c_1 - c_2)) = q n^T (c_1 - c_2)
\]
After local corrections for sharpening terms the assembled sharpening force at the vertex $p_i$ is used in order to correct repulsive force $F_r$:
\[
F_r  \leftarrow F_r - \frac1{2 |F_s|^2} F_s (F_s^T F_r - |F_s^T F_r|)
\]

Attraction force looks like
\[ 
{F}_a(p_i) = -\sum\limits_{e_v: p_i \in \mathop{\rm dual} e_v} \frac12 \left( \frac{L_0}{L} \right)^2 u(c) (C_1 + C_2)^T \frac{\nabla u(c)}{|\nabla u(c)|},
\]

Displacement of Delaunay vertices is done in two steps. The first step is written as
\[
\tilde{p}^0_i = p_i^k + w_r \tau_r F_r + w_s \tau_s F_s, \ 
w_r = \frac1{20}, \ w_s = \frac12
\]
\[
\tau_r = \min (1, \frac{L_0}{ 5 w_r F_r}), \ \tau_s = \min (1, \frac{L_0}{ 5 w_s F_s})
\]
After this displacement we use $M$ iterations with attraction force to project Voronoi edges to boundary
\[
\tilde{p}^{m + 1}_i = \tilde{p}_i^m + \tau_a F_a (\tilde{p}^l_m), \ \tau_a = \frac1{10}
\]
And finally
\[
p^{k+1}_i = \tilde{p}^M_i
\]

\section{Numerical experiments}

We run series of numerical experiments with artificially constructed domains. The complexity of the tests is well represented by the model ``wheel'' shown in Fig.~\ref{fig-domain}. In this model multiple sharp vertices are present on the boundary.

\begin{figure}
\begin{tabular}[c]{cc}
\raisebox{0.35 \textwidth}{(a)} & \includegraphics[width=0.75 \textwidth]{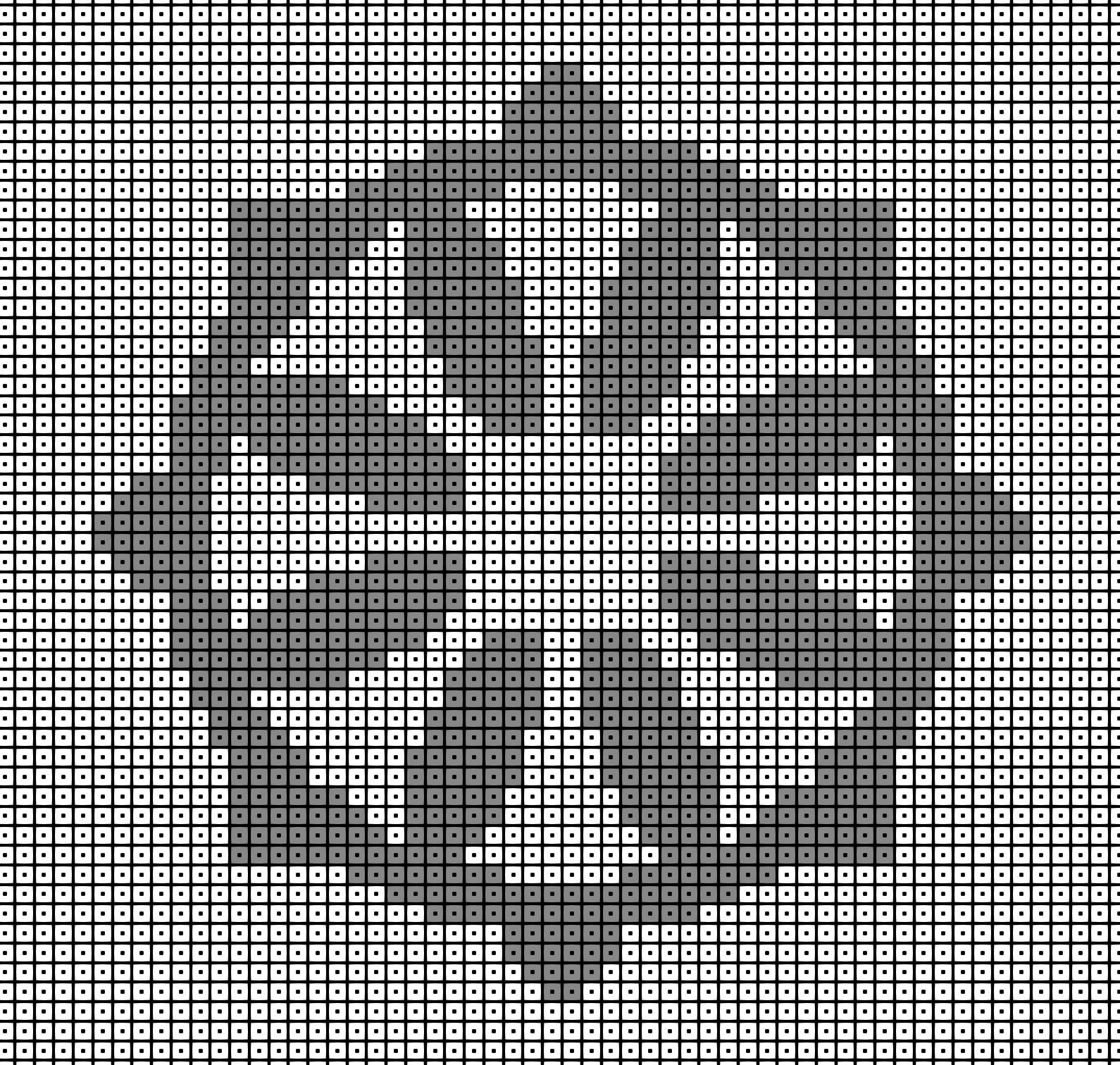} \\
\raisebox{0.35 \textwidth}{(b)} & \includegraphics[width=0.75 \textwidth]{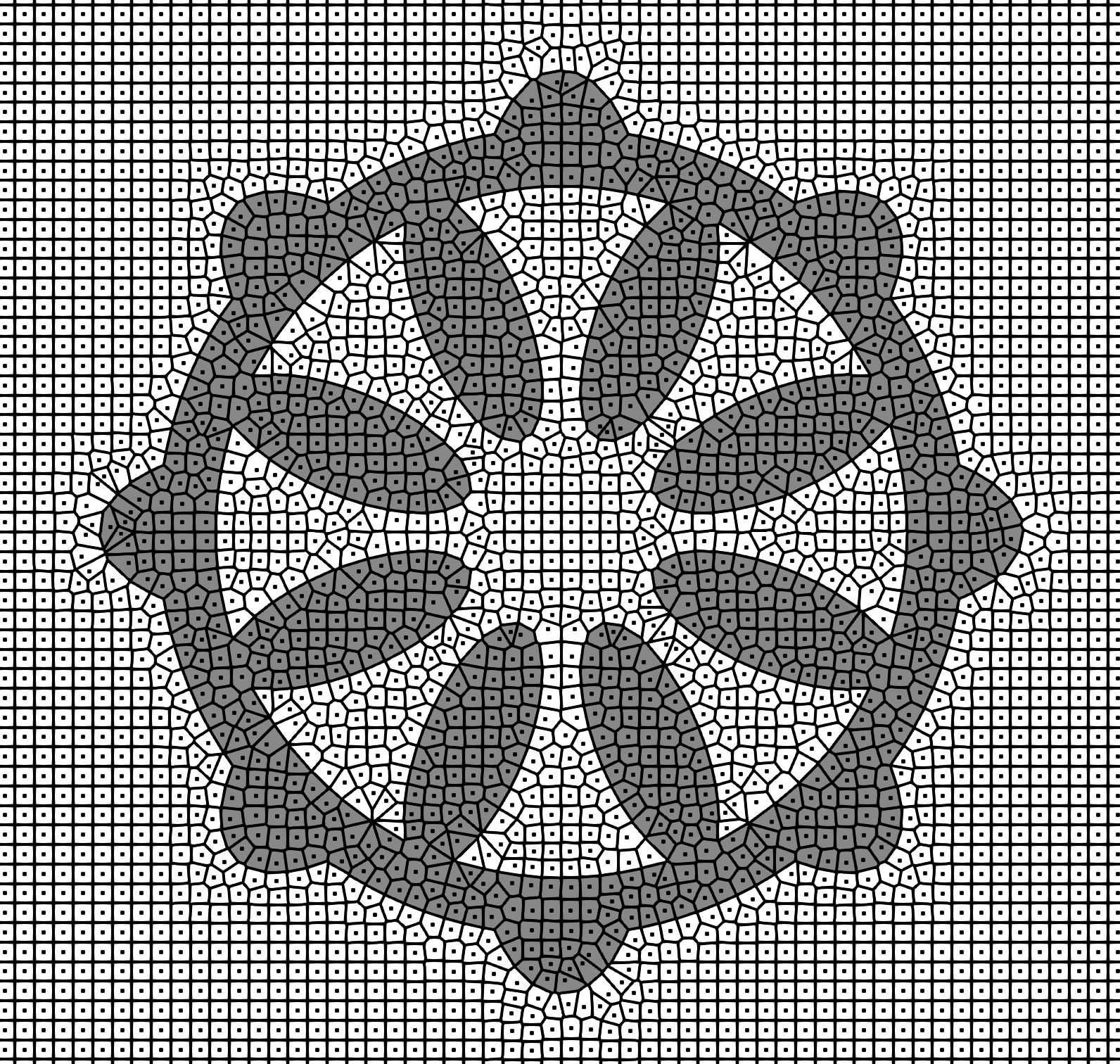}
\end{tabular}
\caption{(a) Initial Voronoi mesh, (b) Voronoi mesh after few iterations.} \label{fig-full0}
\end{figure}

Fig.~\ref{fig-full0} shows initial Cartesian Voronoi mesh and result after few iterations.

\begin{figure}
\centerline{(a) \hspace{0.4\textwidth} (b)}
\includegraphics[width=0.49 \textwidth]{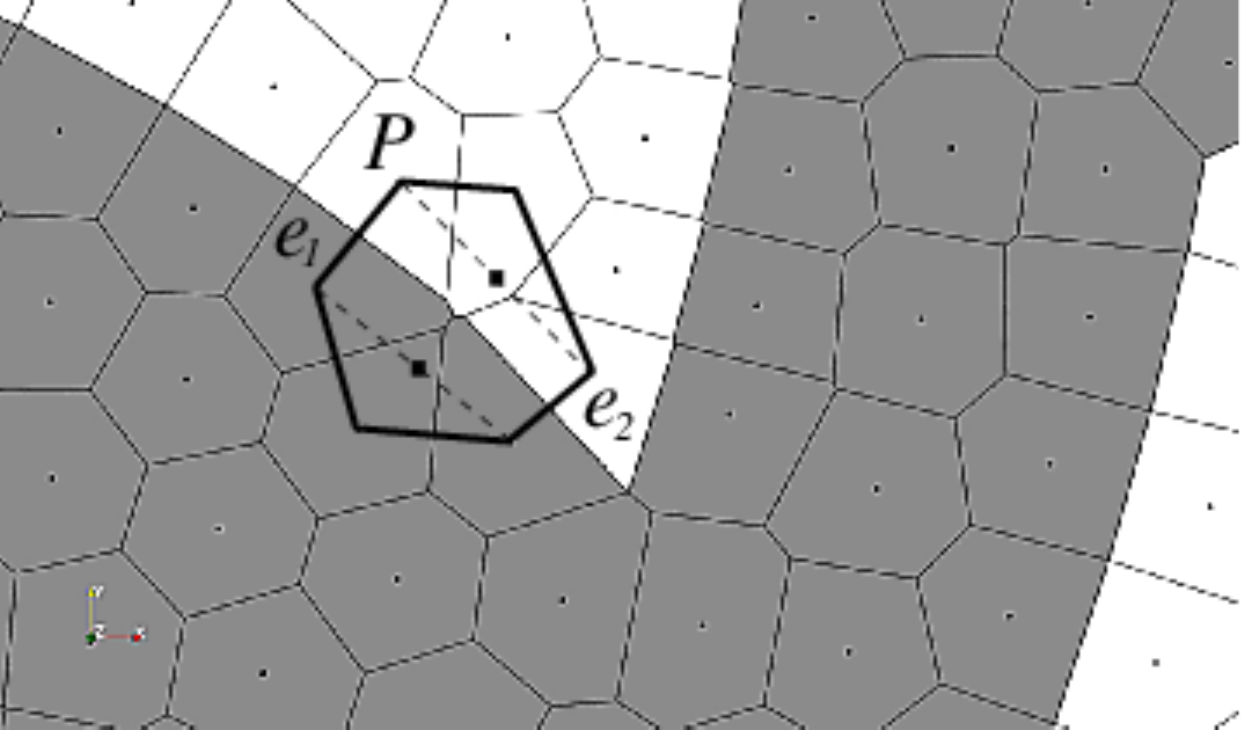} 
\includegraphics[width=0.49 \textwidth]{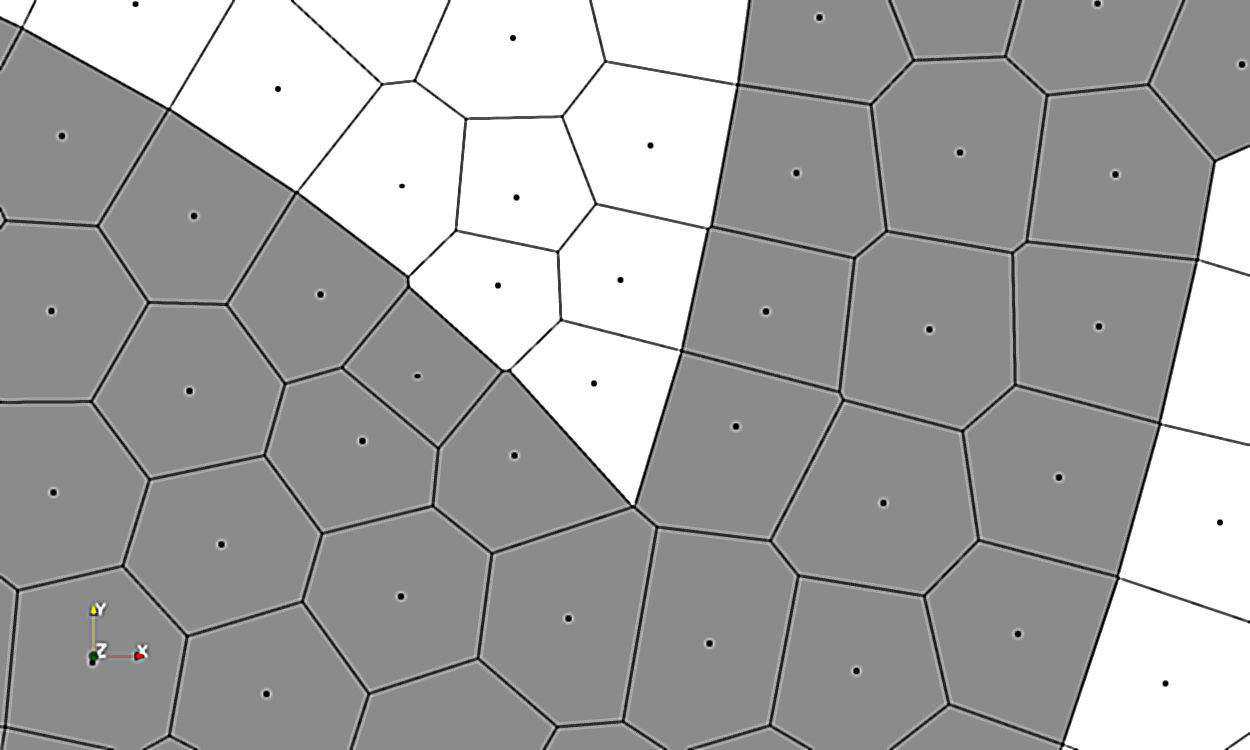}
\caption{(a) Fragment of Voronoi mesh with non-orthogonal edges, (b) correct connectivity is attained by adding new Delaunay vertex.} \label{fig-refine}
\end{figure}

As one can see, algorithm recovers internal boundaries quite fast. However this guess contains approximation defects and layer topological errors when near-boundary Voronoi edges are not orthogonal to boundary. The origin of these errors is simple: Delaunay vertex does not have its mirror across the boundary. Hence most of the topological errors can be eliminated by reasonable Delaunay vertex insertion, as shown in Fig.~\ref{fig-refine}.  We consider polygon $P$ being the closest guess to Delaunay polygon build upon two stable Delaunay edges $e_1$ and $e_2$, crossing the boundary. We build quadrilateral cell upon these two edges and add new vertices at the middle of virtual opposite edges.

Approximate Delaunay hexagon is resolved by inserting two vertices, while approximate Delaunay pentagon is resolved by adding single vertex. In our test cases there was no need to consider more complex polygons.

\begin{figure}
\begin{tabular}[c]{cc}
\raisebox{0.35 \textwidth}{(a)} & \includegraphics[width=0.75 \textwidth]{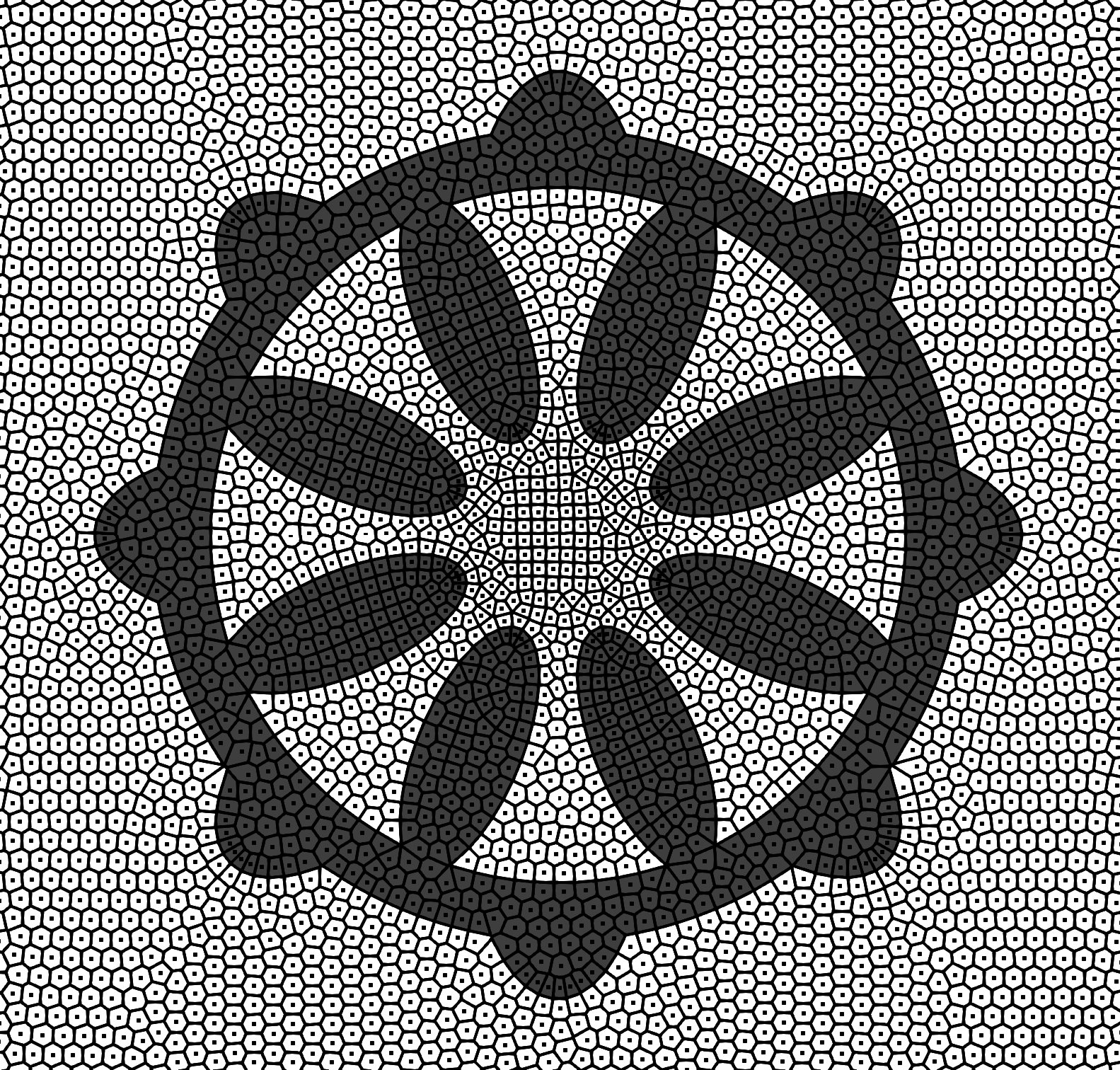} \\
\raisebox{0.35 \textwidth}{(b)} & \includegraphics[width=0.75 \textwidth]{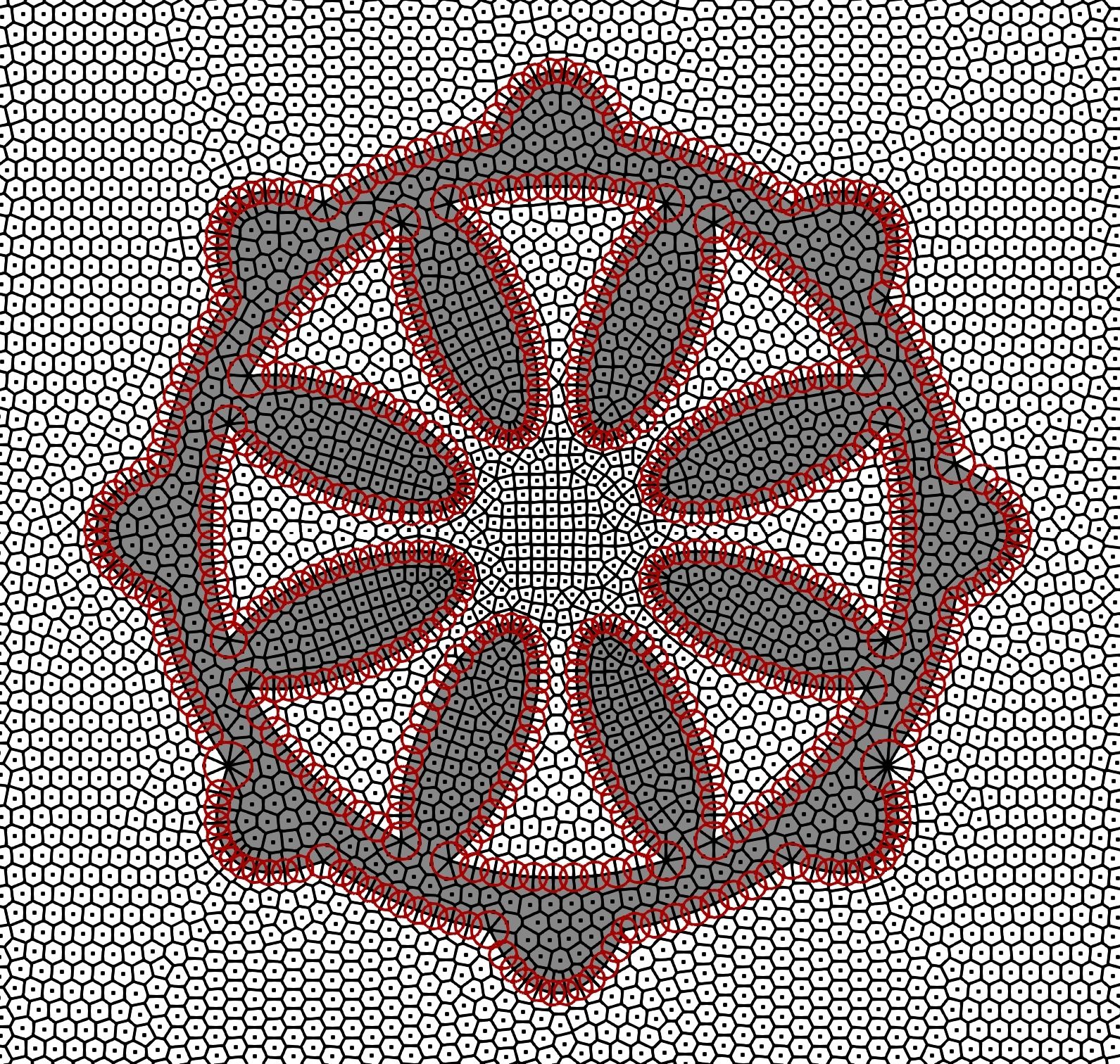}
\end{tabular}
\caption{(a) Stabilized Voronoi mesh, (b) Voronoi mesh after elimination of short edges.} \label{fig-full500}
\end{figure}

Fig.~\ref{fig-full500} shows stabilized Voronoi mesh with fully developed double boundary layer without topological defects. Elimination of small Voronoi edges creates final mesh where internal boundaries are approximated by Voronoi edges and normals to the boundary are approximated by discrete normals.   

\begin{figure}
\includegraphics[width=0.49 \textwidth]{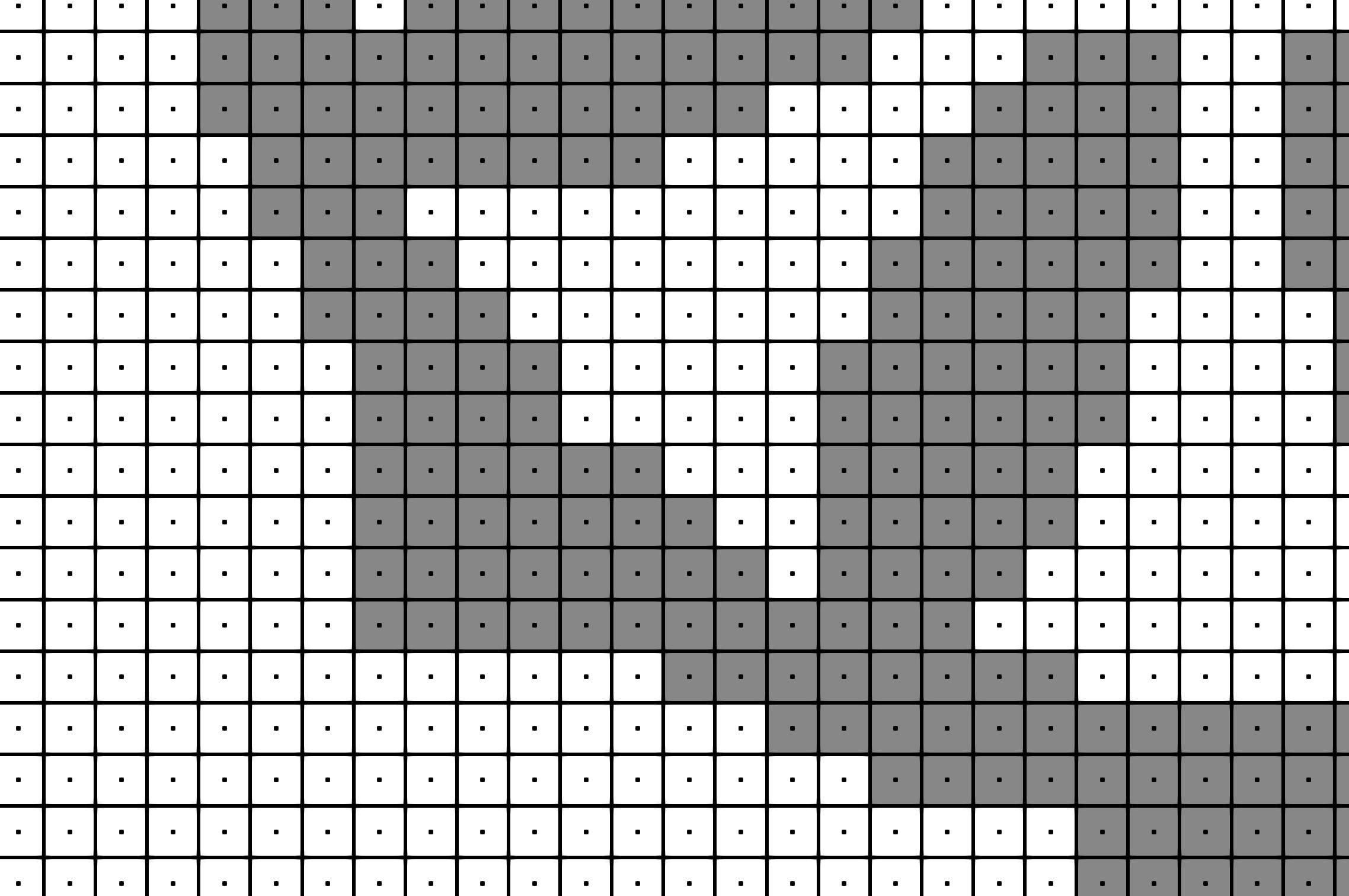}
\includegraphics[width=0.49 \textwidth]{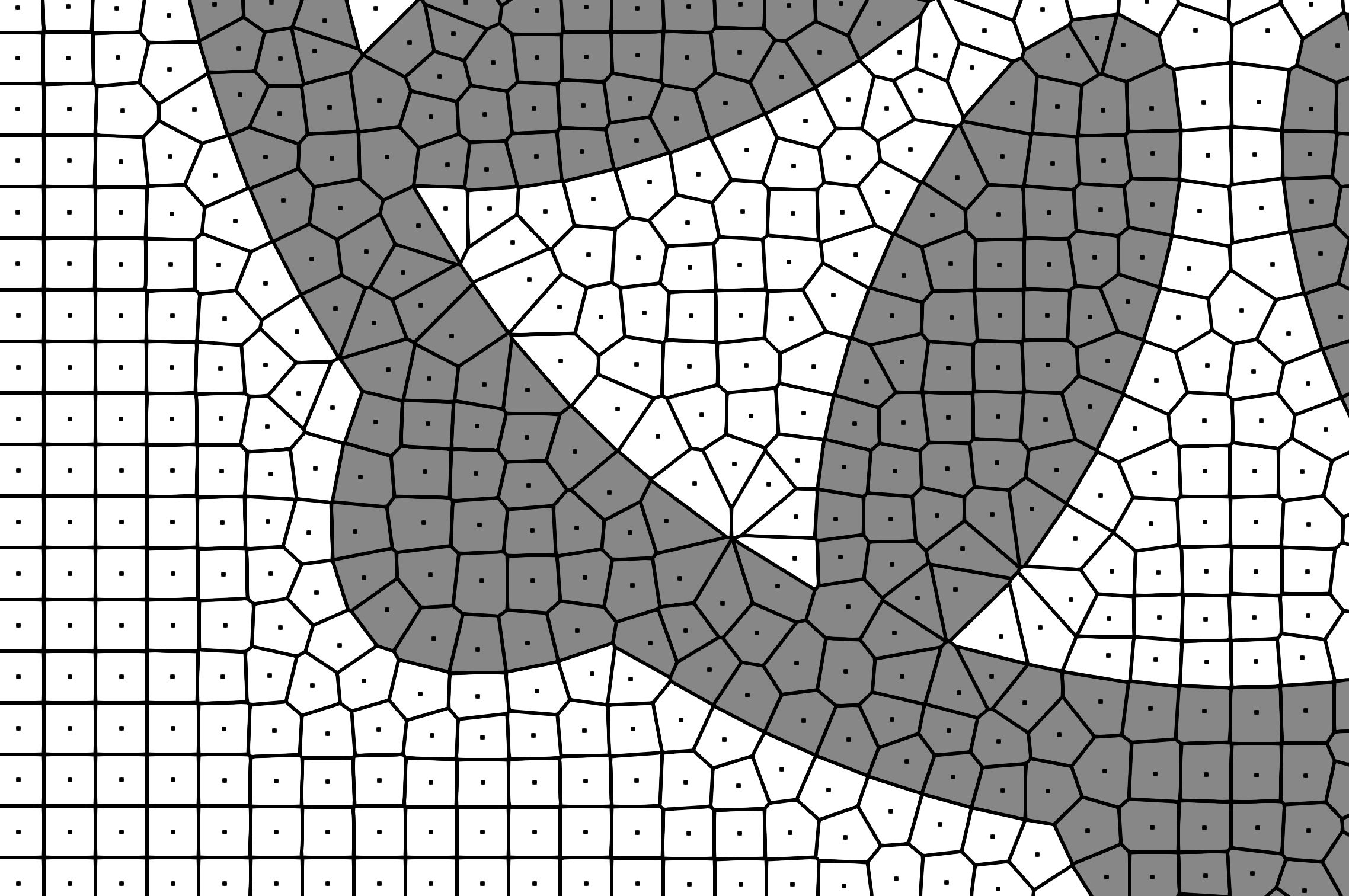}
\caption{Fragment of initial Voronoi mesh and result after few interations.} \label{fig3}
\end{figure}

\begin{figure}
\includegraphics[width=0.49 \textwidth]{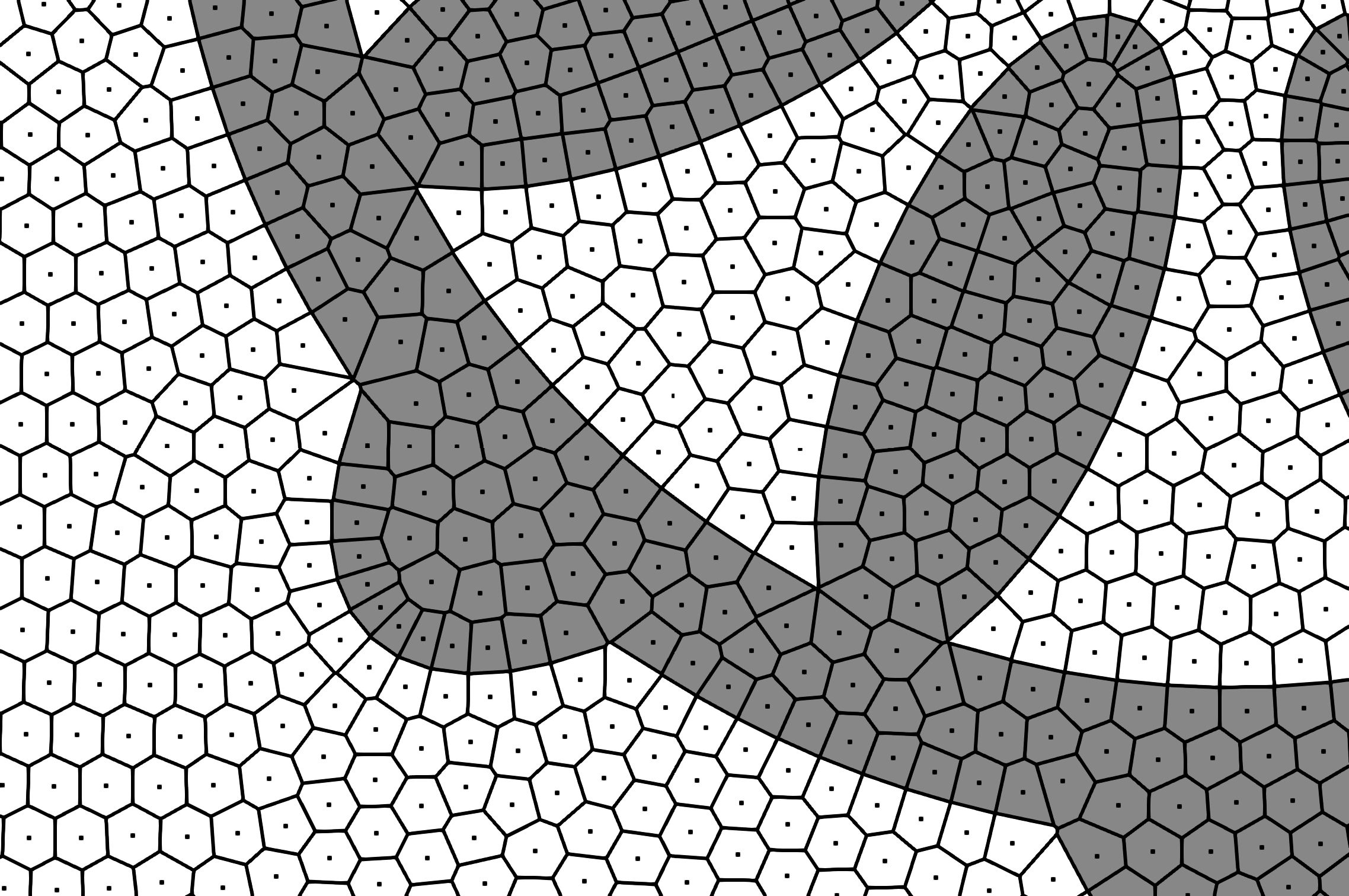}
\includegraphics[width=0.49 \textwidth]{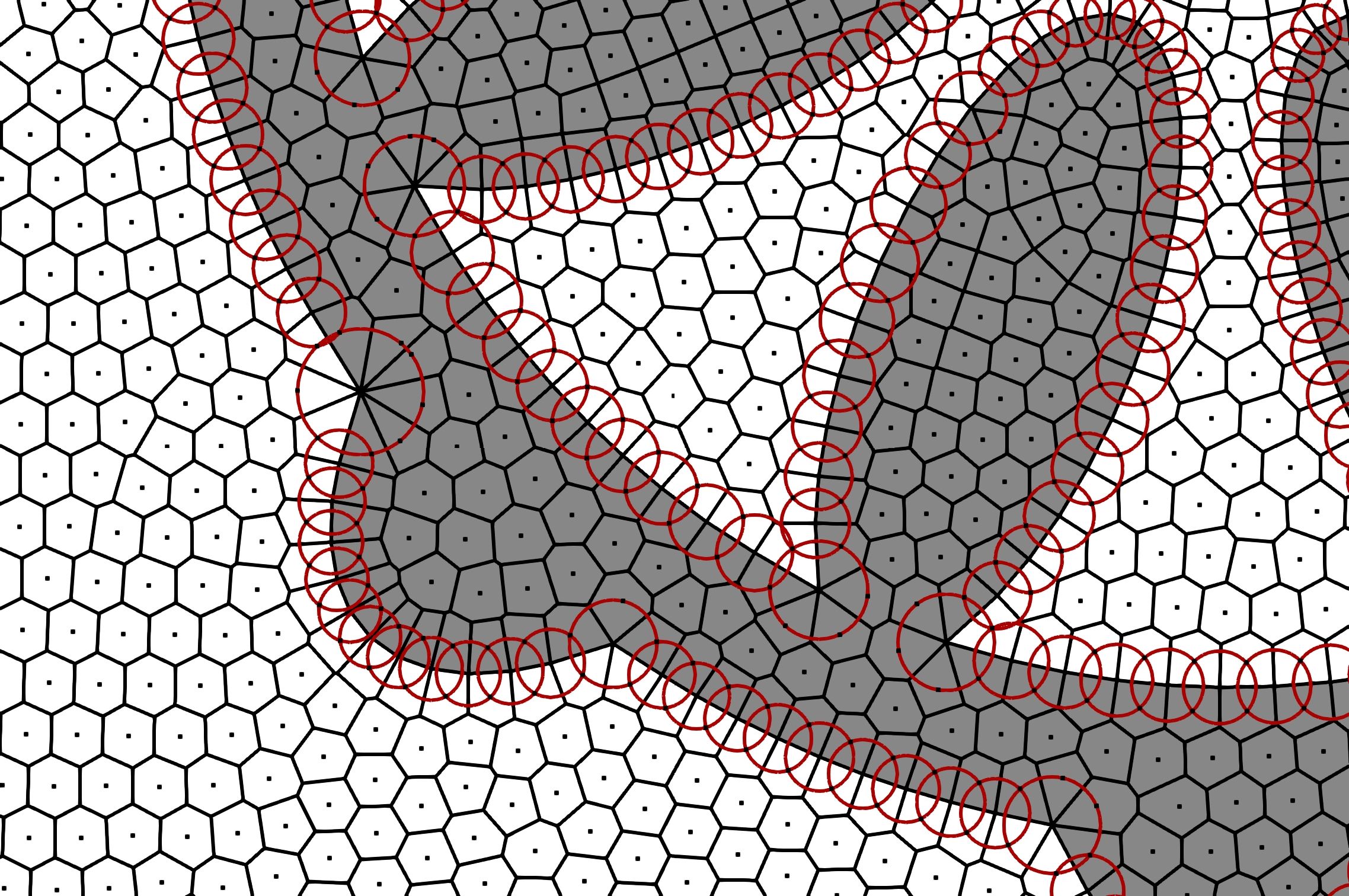}
\caption{Fragment of stabilized Voronoi mesh and result of elimination of short Voronoi  edges.} \label{fig4}
\end{figure}

\begin{figure}
\includegraphics[width=0.49 \textwidth]{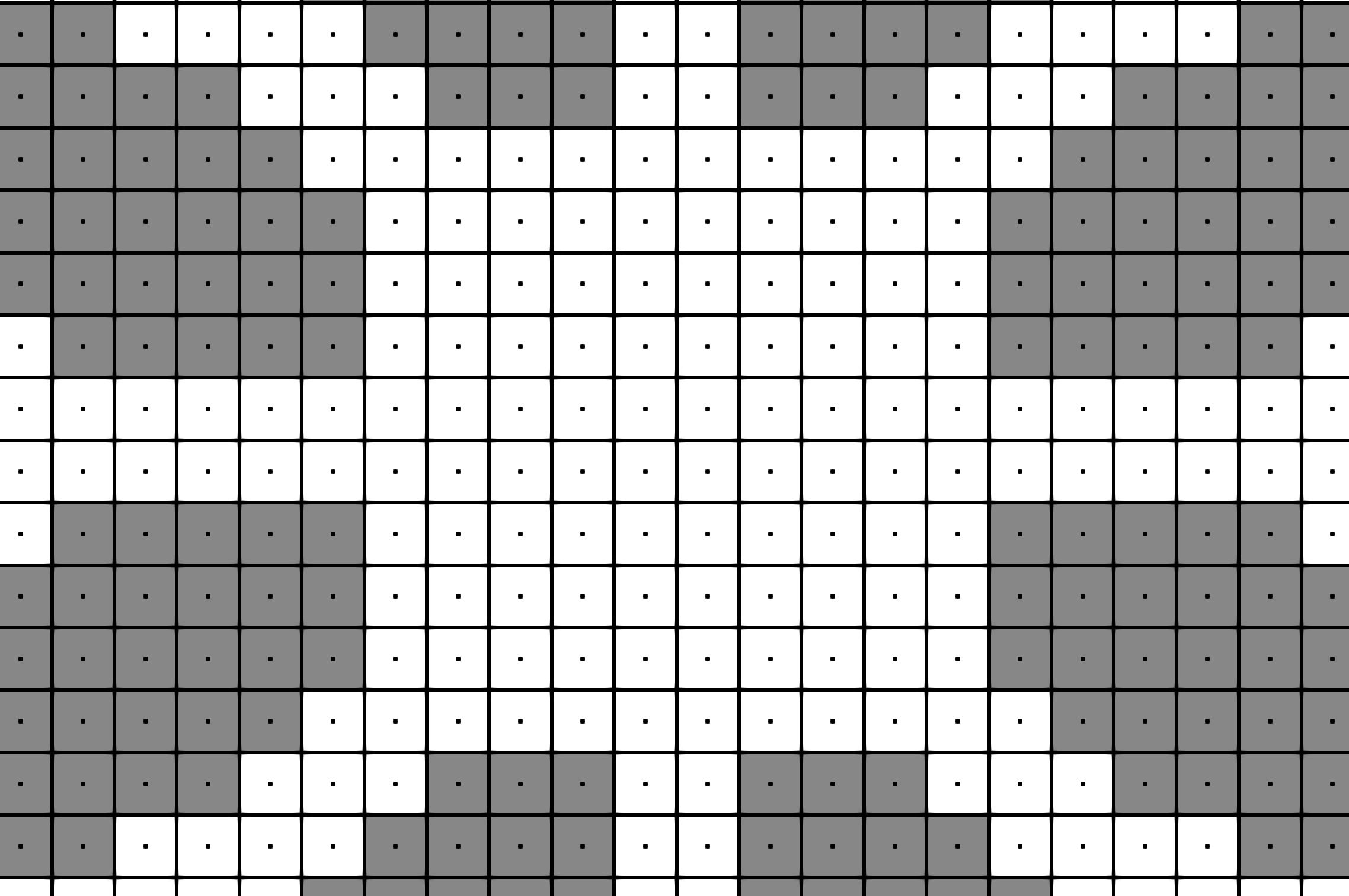}
\includegraphics[width=0.49 \textwidth]{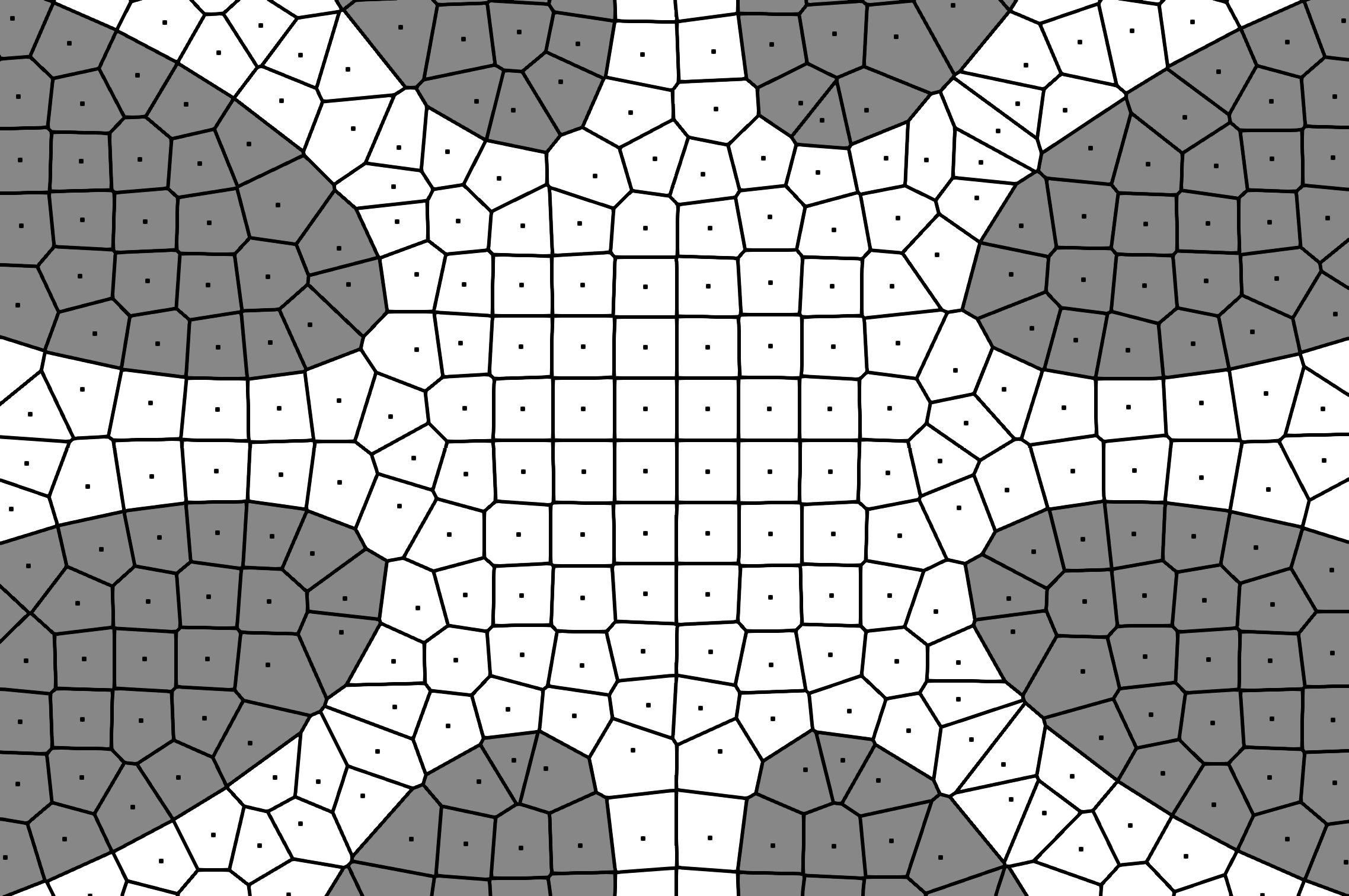}
\caption{Fragment of initial Voronoi mesh and result after few interations.} \label{fig5}
\end{figure}

\begin{figure}
\includegraphics[width=0.49 \textwidth]{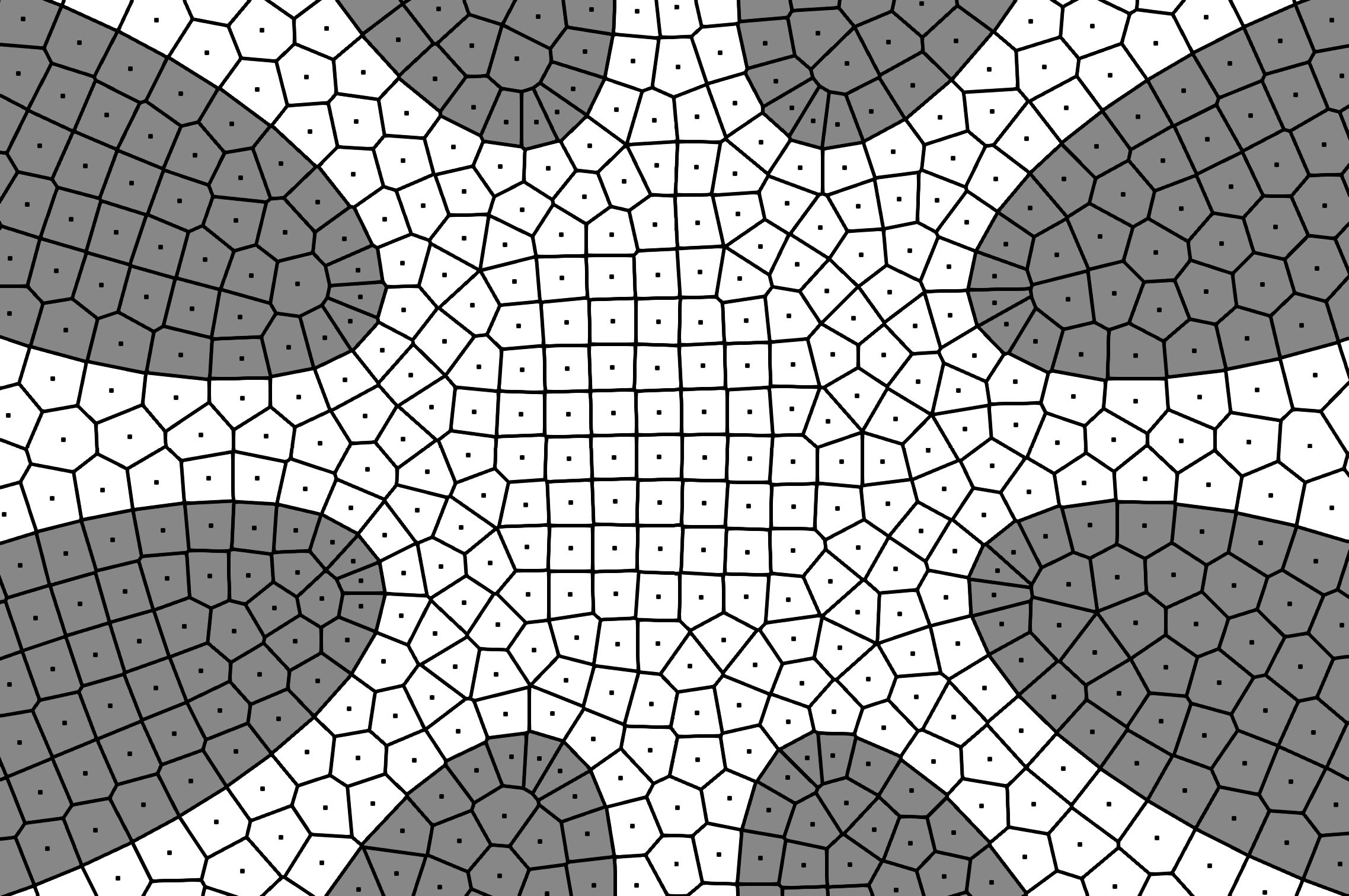}
\includegraphics[width=0.49 \textwidth]{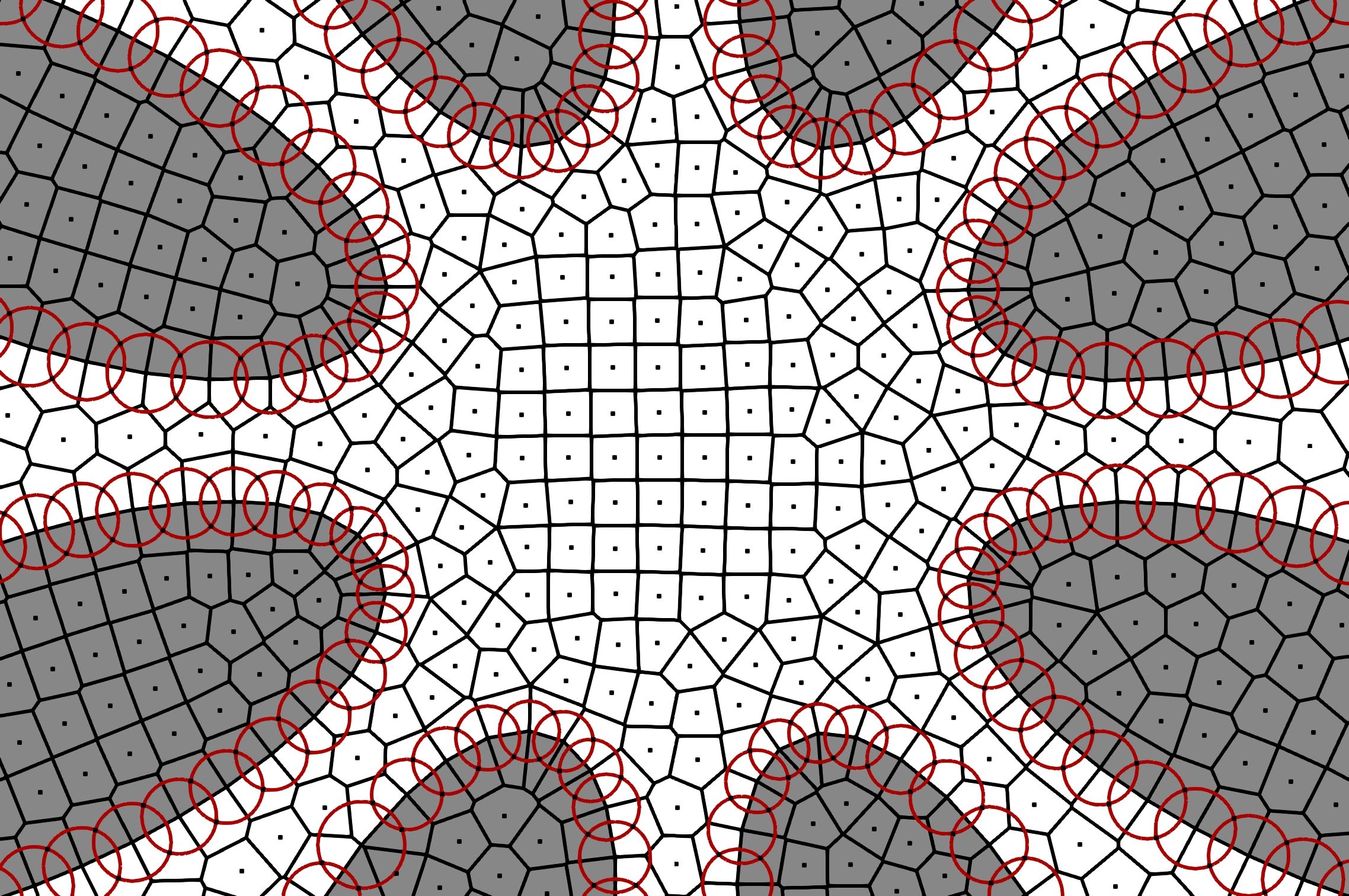}
\caption{Fragment of stabilized Voronoi mesh and result of elimination of short Voronoi  edges.} \label{fig6}
\end{figure}

Figs.~\ref{fig3}-\ref{fig6} illustrate the same step of mesh evolution for two enlarged fragments of the ``wheel'' model.

\begin{figure}
\includegraphics[width=0.49 \textwidth]{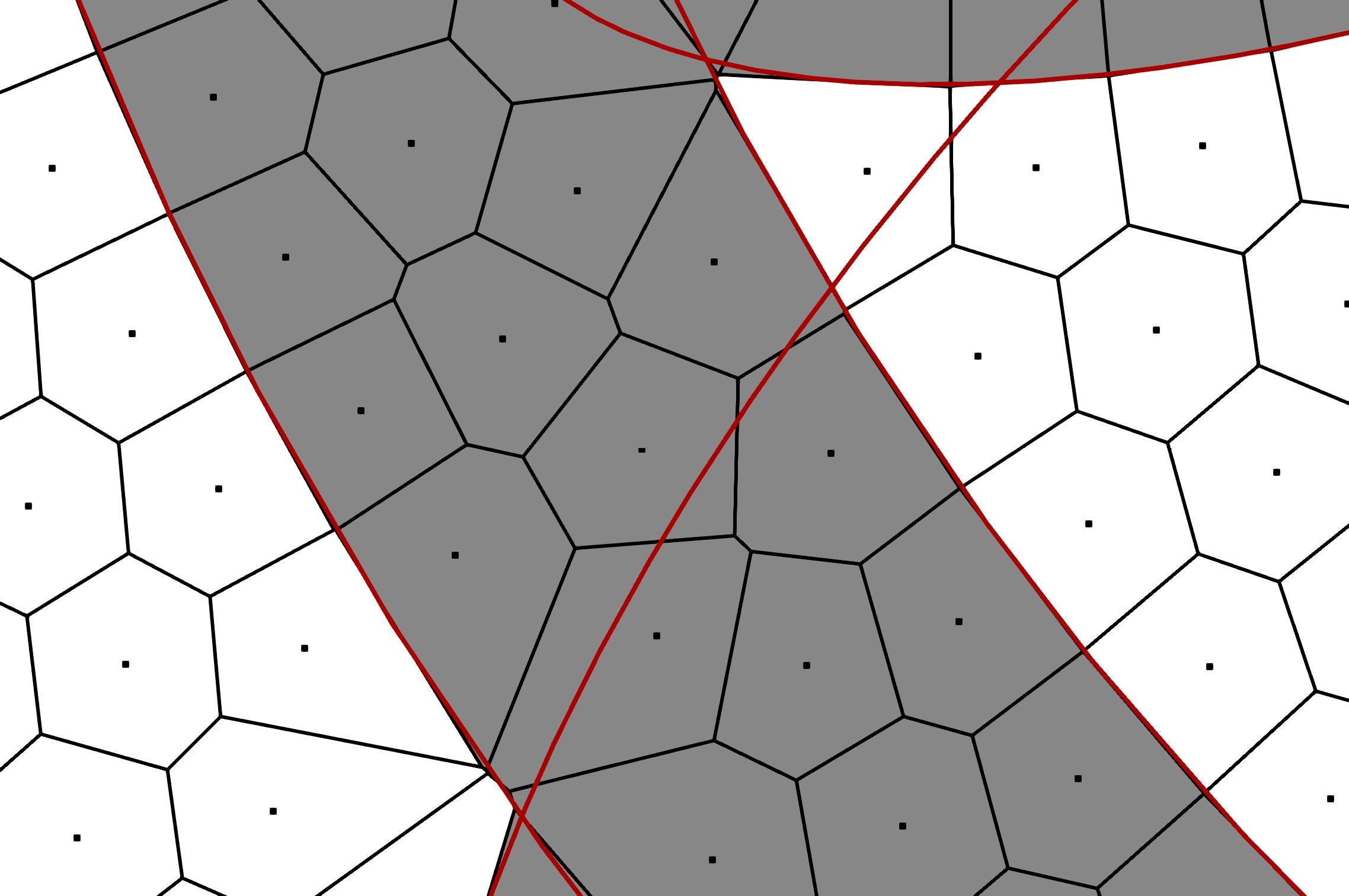}
\includegraphics[width=0.49 \textwidth]{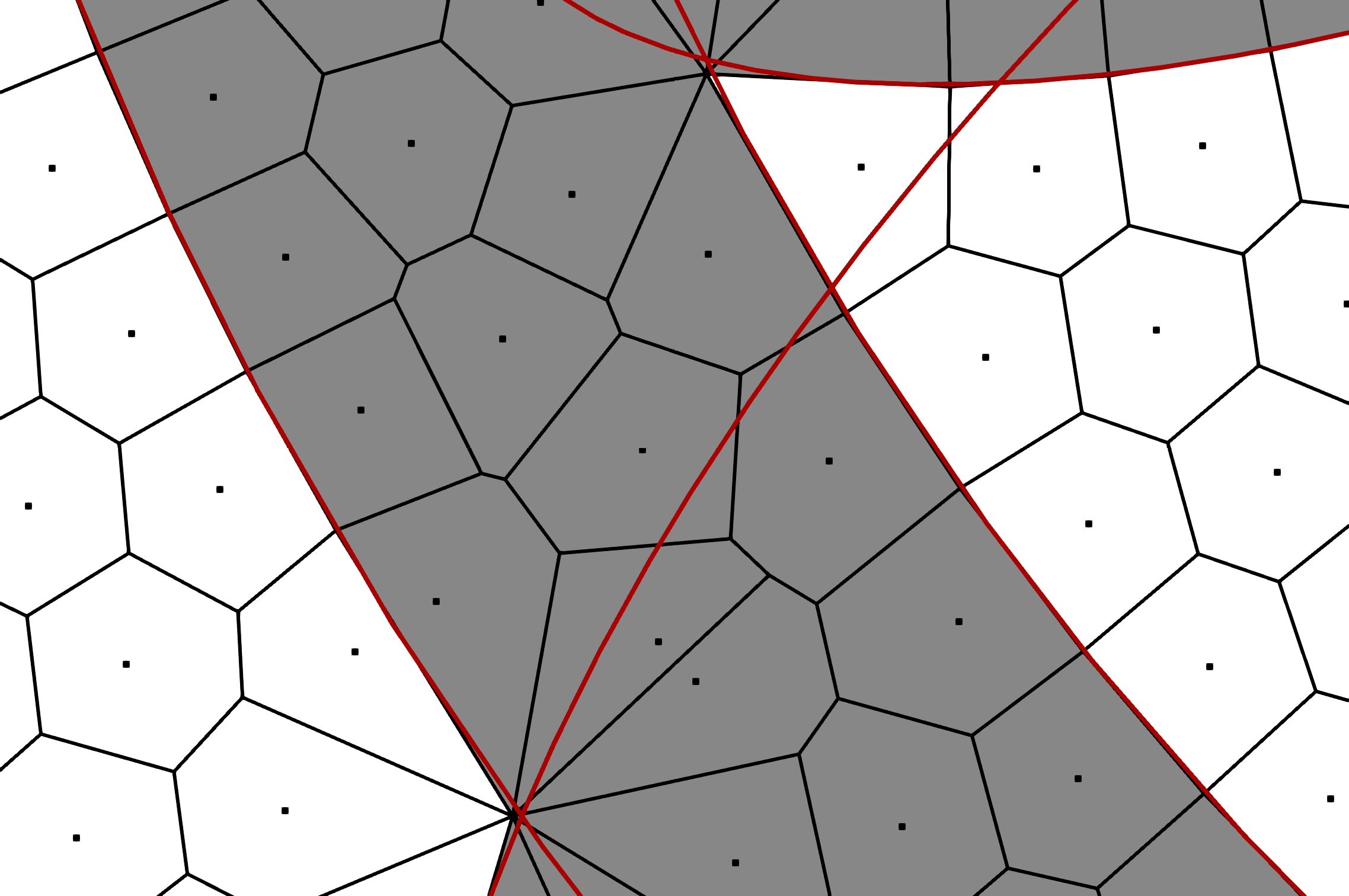}
\caption{Elimination of short Voronoi boundary edges: enlarged view.} \label{fig7}
\end{figure}

\begin{figure}
\includegraphics[width=0.49 \textwidth]{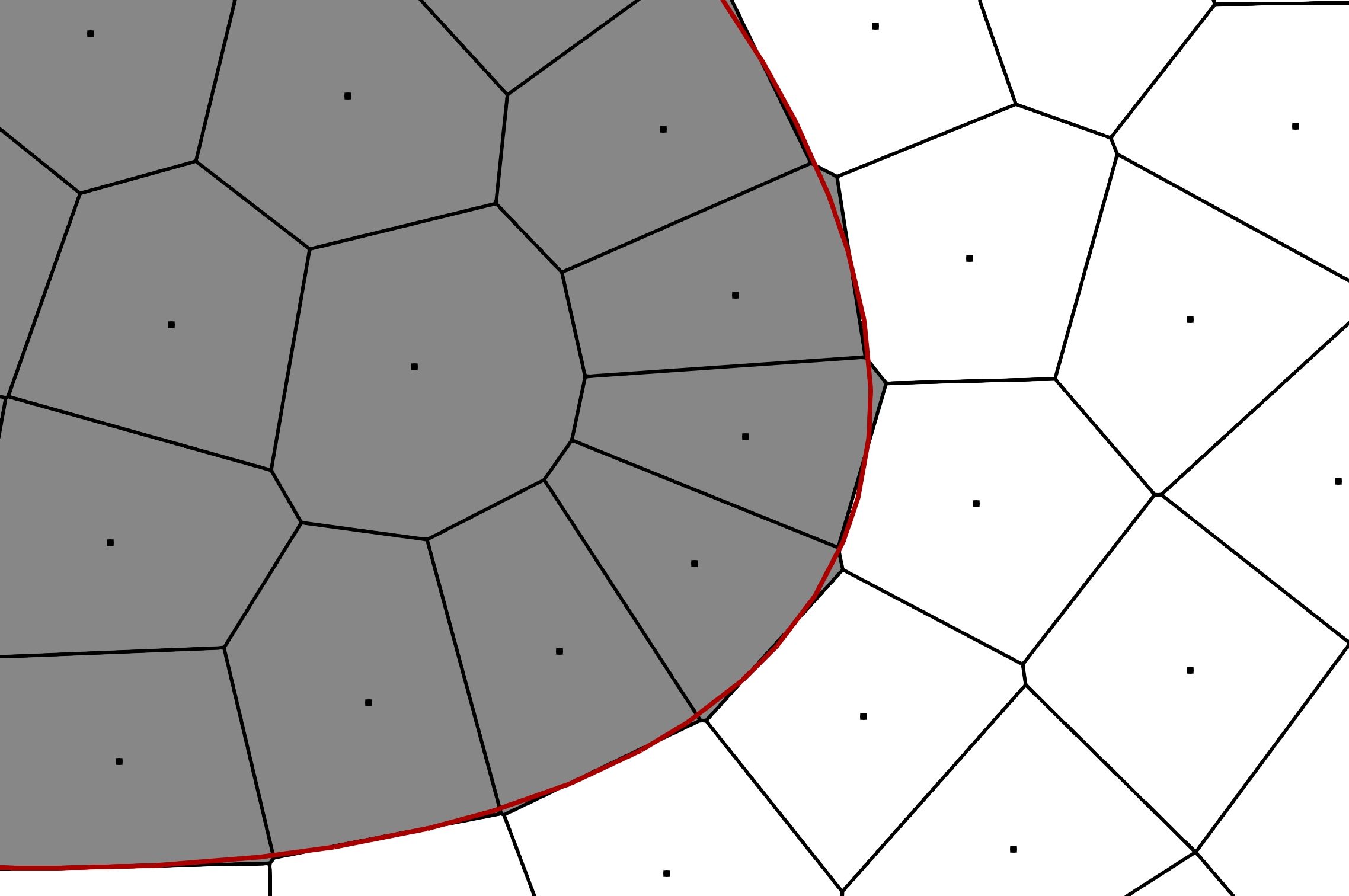}
\includegraphics[width=0.49 \textwidth]{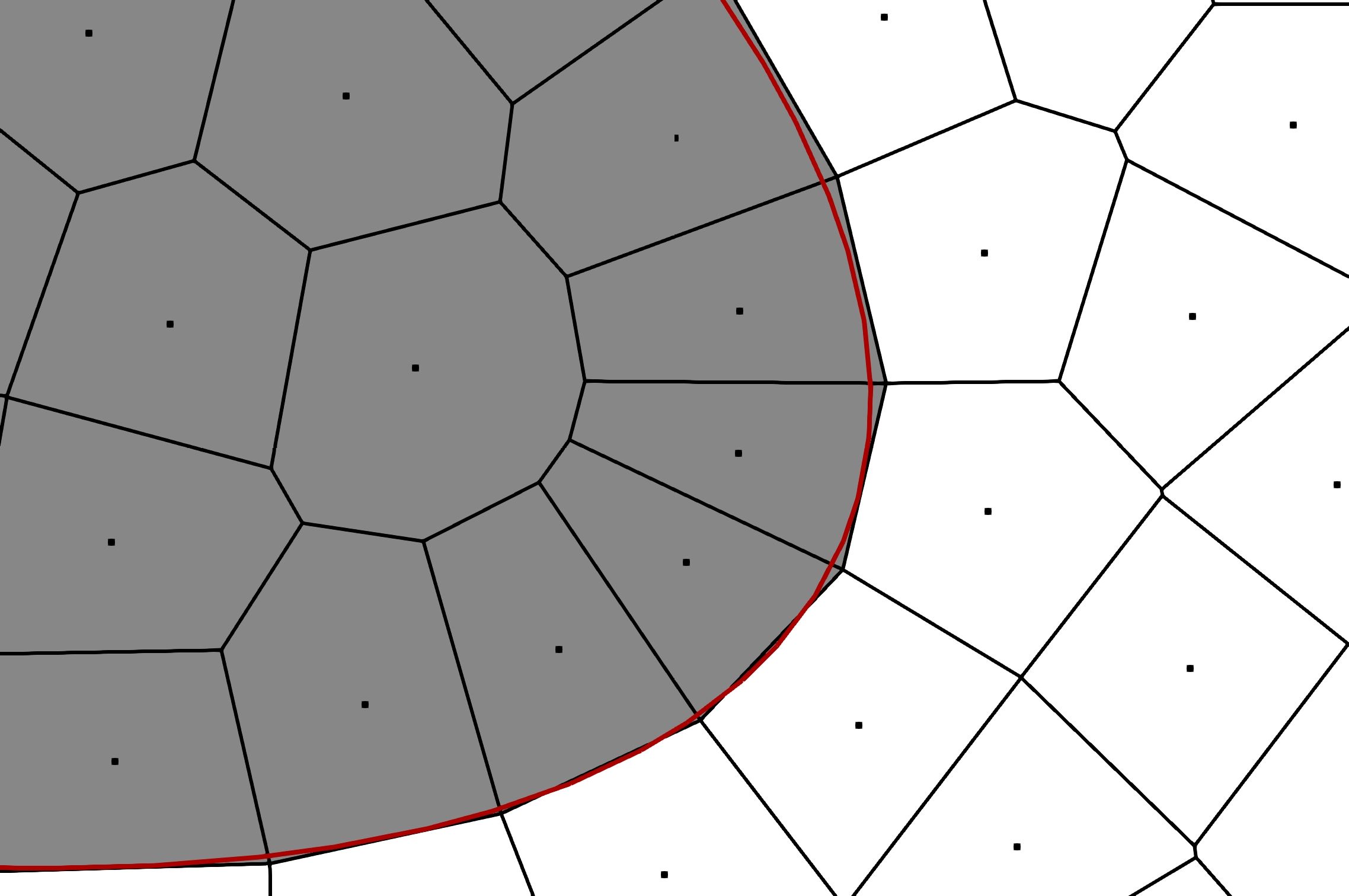}
\caption{Elimination of short Voronoi boundary edges: enlarged view.} \label{fig8}
\end{figure}

Figs.~\ref{fig7},~\ref{fig8} demonstrate that elimination of short Voronoi edges does not lead to deterioration of boundary approximation quality.

\section{Conclusions}

Algorithm for construction of hybrid planar Voronoi meshes demonstrates ability to build orthogonal layers of Voronoi cells near internal boundaries with correct resolution of sharp vertices. It is too early to claim that this algorithm can be used in industrial practice. To this end the following problems should be addressed: current algorithm in some cases creates too large Voronoi cells near sharp corners, case of multimaterial vertices and thin material layers is not yet addressed, multilayered Voronoi layers near boundaries are not considered and most important, generalization to 3d case has to be investigated. These problems are topics of ongoing research.

\end{document}